\documentclass[12pt]{amsart}
\usepackage{amsmath,amssymb,amsthm, amscd}

\newcommand\Z{\mathbb{Z}}
\newcommand\Q{\mathbb{Q}}

\newcommand\calO{{\mathcal{O}}}

\def\p{{\mathfrak{p}}}
\def\q{{\mathfrak{q}}}

\DeclareMathOperator{\GL}{GL}

\newenvironment{Gal}{{\rm Gal\,}}{}

\newenvironment{Aut}{{\rm Aut}}{}
\newenvironment{Frob}{{\rm Frob}}{}

\newtheorem{theorem}{Theorem}[section]
\newtheorem{definition}[theorem]{Definition}
\newtheorem{lemma}[theorem]{Lemma}
\newtheorem{proposition}[theorem]{Proposition}
\newtheorem{proposition-definition}[theorem]{Proposition-Definition}
\newtheorem{corollary}[theorem]{Corollary}
\newtheorem{conjecture}[theorem]{Conjecture}

\theoremstyle{definition}

\theoremstyle{remark}
\newtheorem*{remark}{Remark}
\newtheorem*{example}{Example}

\begin{document}
\title[Iterative Construction]{An Iterative Construction of Irreducible Polynomials Reducible Modulo Every Prime}
\author{Rafe Jones}
\thanks{2010 Mathematics Subject Classification 37P15, 11R09}
\thanks{The author's research was partially supported by NSF grant DMS-0852826.}

\begin{abstract} We give a method of constructing polynomials of arbitrarily large degree irreducible over a global field $F$ but reducible modulo every prime of $F$.  The method consists of finding quadratic $f \in F[x]$ whose iterates have the desired property, and it depends on new criteria ensuring all iterates of $f$ are irreducible.  In particular when $F$ is a number field in which the ideal (2) is not a square, we construct infinitely many families of quadratic $f$ such that every iterate $f^n$ is irreducible over $F$, but $f^n$ is reducible modulo all primes of $F$ for $n \geq 2$.  We also give an example for each $n \geq 2$ of a quadratic $f \in \Z[x]$ whose iterates are all irreducible over $\Q$, whose $(n-1)$st iterate is irreducible modulo some primes, and whose $n$th iterate is reducible modulo all primes.  From the perspective of Galois theory, this suggests that a well-known rigidity phenomenon for linear Galois representations does not exist for Galois representations obtained by polynomial iteration.  Finally, we study the number of primes $\p$ for which a given quadratic $f$ defined over a global field has $f^n$ irreducible modulo $\p$ for all $n \geq 1$. 
\end{abstract}

\maketitle

\section{Introduction} \label{intro}

At the end of the $19$th century, David Hilbert gave examples of irreducible polynomials $f(x) \in \Z[x]$ that are reducible modulo all primes, namely any irreducible member of the family $x^4 + 2ax^2 + b^2$.   In particular, one easily checks that $f(x) = x^4 + 1$ qualifies, since $f(x+1)$ is Eisenstein with respect to 2.  Moreover, $g(x) = x^{2^n}+1$, $n \geq 2$, shares the same properties, since $g(x+1)$ is again Eisenstein and $g(x) = f(x^{2^{n-2}})$ inherits from $f$ a non-trivial factorization modulo any $p$.  In this paper, we give a generalization of this construction, one that yields infinitely many infinite families of irreducible polynomials that are reducible modulo all primes.  Specifically, we give criteria that ensure a quadratic polynomial $f(x) \in \Z[x]$ has its $n$th iterate irreducible over $\Q$ but reducible modulo all primes.  The construction works over most global fields; see Corollary \ref{numfieldcor} and Theorem \ref{funcfield} for exact statements.  Our approach is based on new results dealing with the irreducibility of iterates of quadratic polynomials; see Theorem \ref{numfieldintro}.  For simplicity, we state in Theorems \ref{qintro} and \ref{ratffintro} our results over $\Q$ and $k(t)$, where $k$ is a finite field of odd characteristic.  We denote by $f^n$ the $n$th iterate of a polynomial $f$, and by $\overline{f}$ the coefficient-wise reduction of $f$ modulo a prime.  

\begin{theorem} \label{qintro}
Let $n \geq 2$ and let $f(x) = (x- \gamma)^2 + \gamma + m$, where $m \in \Z$ is arbitrary and $\gamma \in \Z$ is chosen as follows.  Let $f_0(x) = x^2 + m$, and let $s \in \Z$ be a square with $s > (f_0^{n-1}(0))^2$ and with $s$ odd if either $m$ is even or $n$ is odd, and $s$ even otherwise.  Put $\gamma = s - f_0^n(0)$.  Then for any $i \geq n$, $f^i$ is irreducible over $\Q$ and $\overline{f^i}$ is reducible for all primes $p \in \Z$.  
\end{theorem}
For instance, $n = 2$, $m = 0$ and $\gamma = 1$ (coming from $s = 1$) satisfy the hypotheses of the theorem, giving that 
$f(x) = (x - 1)^2 + 1$ has all iterates beyond the first irreducible but reducible modulo all primes.  However, $f^i(x) = (x-1)^{2^i} + 1$, and we recover the example given at the beginning of this section.   Note that Theorem \ref{qintro} applies to $f$ that do not have all iterates Eisenstein.  Take $n = 2$, $m = 1$, and $\gamma = 2$ (this comes from choosing $s = 4$).  Then Theorem \ref{qintro} applies to $f(x) = (x - 2)^2 + 3$, though no iterate of $f$ is Eisenstein since the $x^{2^n - 1}$ coefficient of $f^n$ is a power of two and the constant coefficient is either $0$ or $3$ modulo $4$.  

Our results also allow for the construction of ``primitive" examples where $\overline{f^{n-1}}$ is irreducible for some primes.  In Section \ref{primitive}, for any $n \geq 2$, we construct $f \in \Z[x]$ such that all iterates of $f$ are irreducible over $\mathbb{Q}$, $\overline{f^{n-1}}$ is irreducible for some primes, but $\overline{f^i}$ is reducible for all primes, for $i \geq n$.  For instance, in the case $n = 9$, the polynomial 
\begin{equation} \label{introexample}
f(x) = (x - 88255775491812351975604)^2 + 88255775491812351975605
\end{equation} 
has this property, and indeed there are no similar polynomials with $m, \gamma \in \Z$ having smaller absolute value than those in \eqref{introexample} (see p.~\pageref{minarg}).  Such examples have a natural interpretation in terms of Galois theory.  To $f \in \Z[x]$, associate the {\em arboreal Galois representation} $G_f$, given by the action of the group $\Gal(\overline{\Q}/\Q)$ on the extension of $\Q$ obtained by adjoining all preimages of $0$ under any iterate of $f$.  This set $T$ of preimages, when it does not contain a critical point of $f$, has a natural structure of a rooted tree, with the action of $f$ furnishing the connectivity relation.  The $n$th level of $T$ is the set of vertices of distance $n$ from the root, and these are precisely the roots of $f^n(x)$.  The action of $G_f$ preserves these root sets, and thus preserves each level of $T$.   The results of Section \ref{primitive} imply:

\begin{theorem} \label{galois}
Let $G_f \hookrightarrow \Aut(T)$ be the arboreal Galois representation attached to $f \in \Z[x]$.  Then for each $n \geq 2$ there exists a quadratic $f$ such that $G_f$ acts transitively on each level of $T$, contains an element acting as a $2^{n-1}$-cycle on level $n-1$, and contains no element acting as a $2^n$-cycle on level $n$.  
\end{theorem}

In particular, this implies that the action of $G_f$ on the subtree $T_n \subset T$ consisting of the levels up to $n$ is not as large as possible, since $\Aut(T_n)$ contains $2^n$-cycles.  This suggests a contrast to the case of linear $\ell$-adic representations, that is, homomorphisms $\Gal(\overline{\Q}/\Q) \to \GL_d(\Z_\ell)$, where $\Z_\ell$ denotes the $\ell$-adic integers.  In this case the elements of $(\Z/\ell^n \Z)^d$ may be thought of as the $n$th level of the corresponding tree.  But if the image $G  \leq \GL_d(\Z_\ell)$ of $\Gal(\overline{\Q}/\Q)$ maps onto $\GL_d(\Z/\ell^n\Z)$ for certain small $n$, then $G$ must map onto $\GL_d(\Z/\ell^n\Z)$ for all $n$.  See p. \pageref{galoisdisc} for more discussion.  

The broad applicability of Theorem \ref{qintro} stems from the following new criterion ensuring irreducibility of the iterates of a quadratic polynomial over a number field.  
\begin{theorem} \label{numfieldintro}
Let $F$ be a number field with ring of integers $\calO$, and suppose there is a prime $\q \subset \calO$ with 
$v_\q(2)$ odd.  Let $\gamma, m \in \calO$ and $f(x) = (x- \gamma)^2 + \gamma + m$.  If $\gamma  \not\equiv m \bmod{\q}$ and $-(\gamma + m)$ is not a square in $F$, then $f^n(x)$ is irreducible over $F$ for all $n \geq 1$.  
\end{theorem}
Theorem \ref{numfieldintro} applies to any number field in which the ideal $(2)$ is not a square, and in particular to any number field of odd degree over $\Q$.  The more general version of Theorem \ref{qintro}, Corollary \ref{numfieldcor}, also applies to such fields.  

We now turn to $F = k(t)$, where our result is weaker because we have no equivalent of Theorem \ref{numfieldintro}.  

\begin{theorem} \label{ratffintro}
Let $k$ be a finite field of odd characteristic, $F = k(t)$, and $\calO = k[t]$.  
Let $n \geq 3$ and let $f(x) = (x- \gamma)^2 + \gamma + m$, where $m \in \calO$ has odd degree and $\gamma \in \calO$ is chosen as follows.  Let $f_0(x) = x^2 + m$, and take $\gamma = m^{2^{n-1}} - f_0^n(0)$.  Then $f^n$ is irreducible over $F$ and $\overline{f^n}$ is reducible for all primes $\p \subset \calO$.  
\end{theorem}
We give an example and make some comments on the case $n = 2$ in Section \ref{funcfields}.  When $f$ satisfies the hypotheses of Theorem \ref{ratffintro}, $f^n$ has the curious property that it is irreducible over $k(t)$ but for any $c$ in the algebraic closure of $k$, the specialization of $f$ at $t = c$ is reducible over $k(c)$.  

We note that in \cite{Brandl} and \cite{Guralnick} it is shown that polynomials similar to those in Hilbert's example exist in any composite degree.  These papers adopt a Galois-theoretic viewpoint -- one needs to construct a polynomial whose Galois group acts transitively on the polynomial's roots, but contains no full cycles.  They rely on non-constructive theorems from inverse Galois theory.  Here, we shall not explicitly use the Galois-theoretic perspective except in our treatment of Theorem \ref{galois} in Section \ref{primitive}; for more on the Galois theory of iterates of quadratic polynomials, see e.g. \cite{quaddiv, odoniwn}.  

In Section \ref{setup} we give background and basic results on the irreducibility of iterates of a quadratic polynomial.  In Section \ref{numfields} we prove our main results on number fields, including Theorem \ref{qintro} (see Corollary \ref{qcor}) and Theorem \ref{numfieldintro}.  In Section \ref{primitive}, we construct primitive examples with coefficients in $\Z$ and prove Theorem \ref{galois} (see Theorem \ref{primex}).  In Section \ref{funcfields} we turn to function fields, including Theorem \ref{ratffintro} (see Corollary \ref{ratffcor}).  Finally, in Section \ref{numstab} we study the number of primes $\p$ for which a given quadratic $f$ defined over a global field has $\overline{f}^n$ irreducible for all $n \geq 1$.  The answer should depend on the size and arithmetic of the forward orbit of the critical point of $f$.  We prove this holds when the forward orbit of the critical point is finite or has a certain multiplicative dependence (Theorem \ref{stabnum}), and conjecture that it should be true in the remaining case (Conjecture \ref{stabconj}).  We give a heuristic argument in support of the conjecture and examine some examples.

\section{Setup and Basic Results} \label{setup}

Let $F$ be a field of characteristic $\neq 2$,
and let $f \in F[x]$ be a monic, quadratic polynomial.  By completing the square, we may write
\begin{equation} \label{fdef}
f(x) = (x - \gamma)^2 + \gamma + m.
\end{equation}
Note that $\gamma$ is the unique critical point of $f$.  
\begin{definition}
We call $f \in F[x]$ \textit{stable} if $f^n$ is irreducible over $F$ for all $n \geq 1$.   
\end{definition}
Several recent papers have studied various properties of stable $f$ \cite{ali, ayad, ayadcor, settled, danielson, quaddiv, Shparlinski}.  
The following is one of the fundamental results involving stability, and appears in a slightly different form in \cite[Proposition 3]{settled} (see also \cite[Proposition 4.2]{quaddiv}).

\begin{theorem} \label{fund}
Let $f$ be as in \eqref{fdef}, and let $n \geq 1$.  Then $f^n$ is irreducible if none of $-f(\gamma), f^2(\gamma), f^3(\gamma), \ldots, f^{n}(\gamma)$ is a square in $F$.  Moreover, ``if" may be replaced by ``if and only if" provided that for every finite extension $E$ of $F$ the norm homomorphism $N_{E/F} : E^* \to F^*$ induces an injection $E^*/E^{*2} \to F^*/F^{*2}$.
\end{theorem}

We recall a proof: for $n=1$, we have that $f$ is irreducible if and only if $-f(\gamma)$ is not a square in $F$, since $-f(\gamma) = -(\gamma + m)$.  Let $n \geq 2$ and assume inductively that $f^{n-1}$ is irreducible if none of $-f(\gamma), f^2(\gamma), f^3(\gamma), \ldots, f^{n-1}(\gamma)$ is a square in $F$.  Suppose that none of $-f(\gamma), f^2(\gamma), f^3(\gamma), \ldots, f^{n}(\gamma)$ is a square in $F$.  Then we have $f^{n-1}$ irreducible, and hence separable since $\deg(f^{n-1}) = 2^{n-1}$ and ${\rm char} \; F \neq 2$.  Let $\beta$ be a root of $f^n$, and note that $\alpha := f(\beta)$ is a root of $f^{n-1}$.  Clearly $F(\beta) \supseteq F(\alpha)$.  Now $f^n$ is irreducible if and only if 
$[F(\beta) : F] = \deg(f^n) = 2^n$.  However, $[F(\beta) : F] = [F(\beta) : F(\alpha)] [F(\alpha) : F] = 2^{n-1} [F(\beta) : F(\alpha)]$, where the last equality follows since $f^{n-1}$ is irreducible.  Thus $f^n$ is irreducible if and only if $[F(\beta) : F(\alpha)] = 2$, i.e., if and only if $f(x) - \alpha$ is irreducible over $F(\alpha)$.  We remark that this is a special case of Capelli's Lemma \cite[p. 490]{fein}.  But $f(x) - \alpha$ is irreducible over $F(\alpha)$ if and only if $-(\gamma + m - \alpha)$ is not a square in $F(\alpha)$.  One now computes 
\begin{align} \label{norm}
N_{F(\alpha)/F}(-(\gamma + m - \alpha)) & = \prod_{f^{n-1}(\alpha) = 0} -(\gamma + m - \alpha) \\ 
\nonumber & = (-1)^{2^{n-1}} f^{n-1}(\gamma + m) \\
\nonumber & = f^{n}(\gamma).
\end{align} 
By assumption $f^n(\gamma)$ is not a square in $F$, implying that $-(\gamma + m - \alpha)$ is not a square in $F(\alpha)$ and proving the irreducibility of $f^n$.  In the case where $N_{F(\alpha)/F}$ induces an injection $F(\alpha)^*/F(\alpha)^{*2} \to F^*/F^{*2}$, then $f^n(\gamma)$ is a square in $F$ if and only if $-(\gamma + m - \alpha)$ is a square in $F(\alpha)$, i.e., if and only if $f^n$ is irreducible.  This proves the theorem.

We note that in general $f^n$ will be irreducible even if $f^n(\gamma)$ is a square.  Indeed, in the proof of Theorem \ref{fund}, for $n \geq 2$ we may replace the ground field $F$ by $F_1:=F(\sqrt{-\gamma - m})$, the splitting field of $f$ over $F$.  Then over $F_1$ we have 
$$
f^{n-1}(x) = f(f^{n-2}(x)) = \left(f^{n-2}(x) - \gamma + \sqrt{-(\gamma + m)}\right)\left(f^{n-2}(x) - \gamma - 
\sqrt{-(\gamma + m)}\right).
$$
The two polynomials in the last expression are irreducible over $F_1$ because $f^{n-1}$ is irreducible over $F$, implying that $[F(\alpha) : F_1] = 2^{n-2}$.   Hence \eqref{norm} becomes
\begin{align*}
N_{F(\alpha)/F_1}(-(\gamma + m - \alpha)) & =  (-1)^{2^{n-2}} \left(f^{n-2}(\gamma + m) - \gamma \pm  \sqrt{-(\gamma + m)}\right) \\ 
& = (-1)^{2^{n-2}} \left(f^{n-1}(\gamma) - \gamma \pm  \sqrt{-(\gamma + m)}\right)
\end{align*}
To ease notation, set $\delta = \sqrt{-(\gamma + m)}$, and assume $n \geq 3$.  We now have that $N_{F(\alpha)/F_1}(-(\gamma + m - \alpha))$ is a square in $F_1$ if and only if there are $a, b \in F$ with $(a + b\delta)^2 = f^{n-1}(\gamma) - \gamma \pm \delta.$  This gives $a^2 - b^2(\gamma + m) = f^{n-1}(\gamma) - \gamma$ and $2ab = \pm 1$.  A straightforward computation shows this happens if and only if one of 
\begin{equation} \label{firstnorm}
\frac{1}{2}\left(f^{n-1}(\gamma) - \gamma \pm \sqrt{f^n(\gamma)} \right)
\end{equation}
is a square in $F$.  When $n = 2$ there is an extra minus sign and the elements in question become 
$(-f(\gamma) + \gamma \pm \sqrt{f^2(\gamma)})/2$.  The point of this computation is that the elements in \eqref{firstnorm} may well fail to be squares in $F$ even if $f^n(\gamma)$ is a square.  This observation lies behind our main results, since $f^n(\gamma)$ being a square ensures reducibility of $f^n$ modulo all primes for which $\overline{\gamma}$ and $\overline{m}$ are defined (see Theorem \ref{deddom}).  Because it will be useful to us in the sequel, we state as a theorem:
\begin{theorem} \label{altfund}
Let $f(x) = (x - \gamma)^2 + \gamma + m$ for $\gamma, m \in F$, and let $n \geq 2$.  Then $f^n$ is irreducible if none of $$-f(\gamma), \frac{-f(\gamma) + \gamma \pm \sqrt{f^2(\gamma)}}{2}, \frac{f^2(\gamma) -\gamma \pm \sqrt{f^3(\gamma)}}{2}, \ldots, \frac{f^{n-1}(\gamma) -\gamma \pm \sqrt{f^n(\gamma)}}{2}$$ is a square in $F$.
\end{theorem}

\begin{remark}
The expressions $f^n(\gamma) - \gamma$ are independent of $\gamma$.  Indeed, if we set $f_0(x) = x^2 + m$, then it is easy to see that 
\begin{equation} \label{fzero}
f^n(\gamma) - \gamma = f_0^n(0).  
\end{equation}
\end{remark}

We turn our attention now to Dedekind domains.  The next proposition illustrates the kind of stability result made possible by Theorems \ref{fund} and \ref{altfund}.  It is a mild generalization for quadratic polynomials of a result of Odoni \cite[Lemma 2.2]{odonigalit}, where it is shown that Eisenstein polynomials are stable.  In Theorem \ref{numfield} we give a stronger result in the case where $\calO$ is the ring of integers in a number field.  
\begin{proposition} \label{deddomstab}
Let $\calO$ be a Dedekind domain with field of fractions $F$, let $\gamma, m \in F$, and suppose that there is a prime $\p \subset \calO$ with $v_\p(m)$ positive and odd and $v_\p(\gamma) > v_\p(m)$.  Then $f(x) = (x - \gamma)^2 + \gamma + m$ is stable.  
\end{proposition}

\begin{proof} We use Theorem \ref{fund}.  Note that by \eqref{fzero}, $f^n(\gamma) = f_0^n(0) + \gamma$ for all $n \geq 1$.  Suppose that $v_\p(m) = c$, which is odd and positive by hypothesis; we claim that $v_\p(f_0^n(0)) = c$ for all $n \geq 1$.  For $n = 1$ the claim is clear since $f_0^1(0) = m$.  If $v_\p(f_0^{n-1}(0)) = c$, then $v_\p(f_0^n(0)) = v_\p(f_0^{n-1}(0)^2 + m) = v_\p(m) = c$, where the middle equality follows because $v_\p(f_0^{n-1}(0))^2 = 2c > c$.  As a side note, one can show similarly that if $v_\p(f_0^n(0)) = e > 0$ for any $n$, then $v_\p(f_0^{nm}(0)) = e$ for all $m \geq 1$, or in the terminology of \cite[p. 524]{quaddiv} the sequence $\{(f_0^n(0)) : n \geq 1\}$ is a rigid divisibility sequence. 

We now have that for all $n \geq 1$, $v_\p(f^n(\gamma)) = v_\p(f_0^n(0) + \gamma) = v_\p(f_0^n(0)) = c$, where the middle equality follows since $v_p(\gamma) > v_\p(m)$.  Hence $f^n(\gamma)$ is not a square in $F$.  
\end{proof}

Suppose now that $\calO$ is a Dedekind domain with field of fractions $F$ and that for each $\p \subset \calO$ the residue field $\calO/\p$ is finite.  We recall some basic algebraic facts regarding the ring $\calO_{(c)} := S^{-1}\calO$, where $S = \{c^n : n \geq 0\}$ for some $c \neq 0$ (note that $S$ is multiplicatively closed).  The prime ideals of $\calO_{(c)}$ are precisely those of the form $\p\calO_{(c)}$, where $\p \subset \calO$ does not contain $c$, or equivalently $\p \nmid (c)$.  Moreover, for any such $\p$ we have 
\begin{equation} \label{genred}
\calO_{(c)}/ \p\calO_{(c)} \cong \calO/\p.  
\end{equation}
Now let $f$ be as in \eqref{fdef}, and fix $c \in \calO$ so that $c \gamma  \in \calO$ and $cm \in \calO$. 
Let $R = \calO_{(c)}$, ensuring that $f$ is defined over $R$ (in fact $f$ may be defined over a smaller ring).  Then for each prime $\p \subset \calO$ with $\p \nmid (c)$, \eqref{genred} gives a 
natural ring homomorphism $R \to \calO/\p, x \mapsto \overline{x}$.   By application to coefficients we thus get a polynomial $\overline{f} \in (\calO/\p)[x]$ and $\overline{f^n} = \overline{f}^n$ follows from homomorphism properties.   

\begin{theorem} \label{deddom}
Suppose that $\calO$ is a Dedekind domain with field of fractions $F$ and finite residue fields.  Let $n \geq 2$, let $s \in F$ be a square, and let $m \in F$ be arbitrary.  Put $f_0(x) = x^2 + m$, let $\gamma = s - f_0^n(0)$, and consider $f(x) = (x - \gamma)^2 + \gamma + m$.   
Then $\overline{f^n}$ is reducible for all primes $\p \subset \calO$ with $\p \nmid (c)$, where $c$ satisfies $cs \in \calO$ and $cm \in \calO$.
\end{theorem}


\begin{proof}
We have that $\gamma$ and $m$ belong to $R := \calO_{(c)}$ because $s, m \in R$ and $f_0^n(0)$ is a polynomial in $m$.  Hence $\overline{\gamma}$ and $\overline{m}$ (and in particular $\overline{f}$) are well-defined for all $\p \nmid (c)$.  

For $\p \subset \calO$, the field $F_\p := \calO/\p$ is finite.  For any $\p \nmid (c)$ with $F_\p$ of characteristic $2$, $\overline{f}$ is reducible and hence so is $\overline{f^n}$.  Otherwise $F_\p$ has odd characteristic, and thus any finite extension $E$ of $F_\p$ satisfies $E^*/E^{*2} \cong \Z/2\Z$.  Because $N_{E/F_\p}$ is surjective, the induced map $N_{E/F_\p} : E^*/E^{*2} \to F_\p^*/F_\p^{*2}$ is too, and hence is also injective.  
For $\p \nmid (c)$, we may now write $\overline{f}(x) = (x - \overline{\gamma})^2 + \overline{\gamma} + \overline{m}$ and apply Theorem \ref{fund}.  Using \eqref{fzero} we then have 
$$\overline{f}^n(\overline{\gamma}) = \overline{\gamma + (f^n(\gamma) - \gamma}) = \overline{\gamma + f_0^n(0)} = \overline{s}.$$
By Theorem \ref{fund}, $\overline{f^n}$ is reducible.  
\end{proof}


\section{Results for number fields} \label{numfields}

We now prove Theorem \ref{numfieldintro}, a criterion for stability for certain quadratic polynomials over a number field.  We restate it here.  Denote by $v_\q$ the $\q$-adic valuation for a prime $\q$ of $\calO$.  
\begin{theorem} \label{numfield}
Let $F$ be a number field with ring of integers $\calO$, and suppose there is a prime $\q \subset \calO$ with 
$v_\q(2)$ odd.  Let $\gamma, m \in \calO$ and $f(x) = (x- \gamma)^2 + \gamma + m$.  If $\gamma  \not\equiv m \bmod{\q}$ and $-(\gamma + m)$ is not a square in $F$, then $f$ is stable.  
\end{theorem}

\begin{remark}
The condition on the existence of $\q$ is satisfied provided that the ideal $(2)$ is not the square of another ideal in $\mathcal{O}$. In particular, this must happen when $[F : \Q]$ is odd.  
\end{remark}

\begin{proof}[Proof of Theorem 3.1]  By Theorem \ref{altfund} it suffices to show that $-f(\gamma)$ and all elements of the form 
\begin{equation} \label{crit}
\frac{1}{2} \left(\pm (f^{i-1}(\gamma) - \gamma) \pm \sqrt{f^i(\gamma)}\right), \qquad i \geq 2
\end{equation} 
are not squares in $F$.   Because $f(\gamma) = \gamma + m$, we have that $-f(\gamma)$ is not a square in $F$ by hypothesis.     
If for given $i \geq 2$, $f^i(\gamma)$ is not a square in $F$, then certainly no element of the form \eqref{crit} for the $i$ in question can be a square in $F$.   If $f^i(\gamma)$ is a square in $F$, then we argue as follows.  Suppose that $\q$ divides $\pm (f^{i-1}(\gamma) - \gamma) \pm \sqrt{f^i(\gamma)}$, so that $\pm (f^{i-1}(\gamma) - \gamma) \equiv \pm \sqrt{f^i(\gamma)} \bmod{\q}$.  Squaring and using \eqref{fzero} then gives $f_0^{i-1}(0)^2 \equiv f^i(\gamma) \bmod{\q}$.  Hence $f_0^{i}(0) - m \equiv f^i(\gamma) \bmod{\q}$, and applying \eqref{fzero} again yields
$$f^{i}(\gamma) - \gamma - m \equiv f^i(\gamma) \bmod{\q}.$$
Because $\calO/\q$ has characteristic two, this implies that $\gamma \equiv m \bmod{\q}$, a contradiction.  


We now have 
$$v_\q \left(\frac{\pm(f^{i-1}(\gamma) -\gamma) \pm \sqrt{f^i(\gamma)}}{2}\right) = v_\q(1/2) = -v_\q(2),$$
and the latter is odd, showing that none of the elements of the form \eqref{crit} is a square in $F$.  
\end{proof}

\begin{corollary} \label{numfieldcor}
Let $F$ be a number field with ring of integers $\calO$, and suppose there is a prime $\q \subset \calO$ with 
$v_\q(2)$ odd.  Let $n \geq 2$, fix $m \in \calO$, let $f_0(x) = x^2 + m$, and choose $s \in \calO$ to be a square such that $s - (f_0^{n-1}(0))^2 \not\equiv 0 \bmod{\q}$ and $-(s - (f_0^{n-1}(0))^2)$ is not a square in $F$.  Then putting $\gamma = s - f_0^n(0)$ and $f(x) = (x- \gamma)^2 + \gamma + m$ we have that for any $i \geq n$, $f^i$ is irreducible over $F$ and $\overline{f^i}$ is reducible for all $\p \subset \calO$.  
\end{corollary}
 
\begin{proof} Note that $\gamma + m = s - f_0^n(0) + m = s - (f_0^{n-1}(0))^2$, and so the hypotheses imply 
that $\gamma + m \not\equiv 0 \bmod{\q}$ and $-(\gamma + m)$ is not a square in $F$.  By Theorem \ref{numfield}, $f$ is stable, and so in particular $f^i$ is irreducible for all $i \geq n$.  On the other hand, since $m,s \in \calO$ we may take $c = 1$ in Theorem \ref{deddom}, showing that $\overline{f^n}$ is reducible for all $\p \subset \calO$.  Then $\overline{f^i} = \overline{f^n} \circ \overline{f^{i-n}}$, which is reducible for all $\p \subset \calO$.  
\end{proof}

\begin{remark} For each $m \in \calO$ it is possible to find infinitely many values of $s$ satisfying the hypotheses of Corollary \ref{numfieldcor}.  Indeed, fix a prime $\mathfrak{r}$ of $\calO$ not dividing $(2)$ or $(f_0^{n-1}(0))$, and let $x \in \mathfrak{r} / \mathfrak{r}^2$.  By the Chinese remainder theorem there exist infinitely many $a \in \calO$ with $a \equiv f_0^{n-1}(0) + x \bmod{\mathfrak{r}^2}$ and $a \not\equiv f_0^{n-1}(0) \bmod{\q}$.  Taking $s = a^2$ satisfies the hypotheses of Corollary \ref{numfieldcor}.  To see why, note that $a + f_0^{n-1}(0) \equiv 2f_0^{n-1}(0) \not\equiv 0 \bmod{\mathfrak{r}}$, and so $\mathfrak{r}$ divides $s - f_0^{n-1}(0)^2$ to only the first power, showing it is not a square in $F$.  Also, $a \not\equiv f_0^{n-1}(0) \bmod{\q}$ implies $a \not\equiv -f_0^{n-1}(0) \bmod{\q}$ since $\q \mid (2)$, and so 
$s - f_0^{n-1}(0)^2 \not\equiv 0 \bmod{\q}$.
\end{remark}

\begin{corollary} \label{qcor}
Fix $n \geq 2$ and $m \in \Z$, and let $s \in \Z$ be a square with $s$ odd if either $m$ is even or $n$ is odd, and $s$ even otherwise.   Let $f_0(x) = x^2 + m$, and suppose that $s > (f_0^{n-1}(0))^2$.  Then putting $\gamma = s - f_0^n(0)$ and $f(x) = (x- \gamma)^2 + \gamma + m$ we have that for any $i \geq n$, $f^i$ is irreducible over $F$ and $\overline{f^i}$ is reducible for all primes $p \in \Z$.  
\end{corollary}

\begin{proof} By Corollary \ref{numfieldcor}, we only need to show that $s - (f_0^{n-1}(0))^2$ is odd and 
$-(s-(f_0^{n-1}(0))^2)$ is not a square in $\Q$.  The latter is immediate from $s > (f_0^{n-1}(0))^2$, while the former follows from the observation that $f_0^{n-1}(0)$ is even if $m$ is even or $n$ is odd, and odd otherwise.  
\end{proof}


For a given $m$, Corollary \ref{qcor} can be used to find infinitely many $\gamma$ such that $f(x)$ is stable but $\overline{f^n}$ is reducible for all primes, for any $n \geq 2$.  Indeed, let $n = 2$ and choose $s$ of parity and size satisfying the hypotheses of Corollary \ref{qcor}.  For instance, when $m = 0$ any odd $s$ will do, though the resulting polynomials $f(x) = (x - s)^2 + s$ have iterates with the closed form $f^n(x) = (x-s)^{2^n} + s$.  For a family whose iterates do not have a closed form, let $m = 1$; then $n = 2$ implies we need to take $s$ even with $s > 1$.  Setting $s = (2a)^2$ with $a \in \Z, a \geq 1$ gives $\gamma = s - f_0^2(0) = 4a^2 - 2$ and this yields the family 
$$f(x) = (x - \gamma)^2 + \gamma + 1= x^2 + (-8a^2 + 4)x + 16a^4 - 12a^2 + 3, \qquad a \geq 1$$ 
any member of which is stable but has $\overline{f^n}$ reducible for all primes, for any $n \geq 2$.  
Many more examples can be found in the next section.  

\section{Primitive examples} \label{primitive}

We can use Corollary \ref{qcor} to generate ``primitive" examples, namely where $f$ is stable, $\overline{f^n}$ is reducible for all primes, and $\overline{f^{n-1}}$ is irreducible for some primes.  
For instance, let $n = 9$ and $m = 1$.  We have $$f_0^{9}(0) = 1947270476915296449559703445493848930452791205
.$$
Set $s = (f_0^8(0) + 1)^2$, which is odd and thus satisfies the hypotheses of Corollary \ref{qcor}.  We then have 
\begin{equation} \label{gam}
\gamma = s - f_0^9(0) = 88255775491812351975604,
\end{equation}
and thus by Corollary \ref{qcor}, the 9th iterate of the polynomial 
$$f(x) = (x - 88255775491812351975604)^2 + 88255775491812351975605$$ is irreducible over $\Q$ but reducible modulo all primes $p$.   By Theorem \ref{fund}, $\overline{f^8}$ is irreducible for any $p$ such that none of $-f(\gamma), f^2(\gamma),f^3(\gamma), \ldots, f^8(\gamma)$ is a square modulo $p$.
Using a computer, one verifies the following condition:
\begin{enumerate}
\item[(*)] For each $1 \leq i \leq 8$ there is an odd prime $r_i$ dividing $f^i(\gamma)$ to odd multiplicity, and, when $i \geq 2$, not dividing $f^k(\gamma)$ for $1 \leq k < i$.
\end{enumerate}
Using quadratic reciprocity and the Chinese remainder theorem one can find $p$ such that $r_i$ is not a square modulo $p$ but each of $-f(\gamma)/r_1, f^2(\gamma)/r_2, \ldots, f^8(\gamma)/r_8$ is a square modulo $p$.  Then $\overline{f^8}$ is irreducible for this $p$.  Indeed, condition (*) implies that the numbers $-f(\gamma), f^2(\gamma),f^3(\gamma), \ldots, f^8(\gamma)$ are linearly independent in the $\Z/2\Z$-vector space $\Q^*/\Q^{*2}$, and using Kummer theory and the Chebotarev density theorem one obtains that the density of primes $p$ for which $\overline{f^8}$ is irreducible is $2^{-8}$.  

\label{minarg} Moreover, if $f$ is any quadratic polynomial with $m, \gamma \in \Z$, $\overline{f^8}$ irreducible for some primes, and $f^9(\gamma)$ a square (whence 
$\overline{f^9}$ is reducible for all primes), then $|m| \geq 1$ and $|\gamma|$ is at least the value given in \eqref{gam}.  To see why, note that we cannot have $m \in \{-2,-1,0\}$, for otherwise the set $\{f^n(\gamma) : n \geq 1\}$ is finite, and by the proof of Theorem \ref{stabnum} it follows that either $f^2$ is reducible modulo all primes or there is a prime with $\overline{f^n}$ irreducible for all $n$.  Now $f_0^2(0) = m(m+1)$, which is at least $2m$ if $m > 0$, and 
at least $|m|(|m|-1)$ otherwise.  Hence if $m \not\in \{-2,-1,0\}$, then $f_0^2(0) \geq 2|m|$, and it is easy to see that this gives $f_0^n(0) > 2|m|$ for $n \geq 3$.  
 
We observe now that $s = (f_0^8(0) + c)^2$ for some $c \in \Z$, implying that $\gamma = s - f_0^9(0) = 2cf_0^8(0) + c^2 - m$.  Fixing $m$, we see that $\gamma$ is quadratic in $c$, and hence the integer $c$-values that minimize $|\gamma|$ must be the nearest integers to 
\begin{equation} \label{roots}
c = -f_0^8(0) \pm \sqrt{f_0^8(0)^2 + m},
\end{equation}
which are the zeroes of $\gamma$.  It is straightforward to verify that if $y > 2|m|$ and $|m| \geq 1$, then 
$$y - 1/2 < \sqrt{y^2 + m} < y + 1/2,$$ 
and thus the integers nearest the roots in \eqref{roots} are $0, \pm 1, -2f_0^8(0),$ and $-2f_0^8(0) \pm 1$.  
We cannot have $c = 0$ or $c = -2f_0^8(0)$, for then $\gamma = -m$, and $f(x) = (x - \gamma)^2$ is already reducible, and hence so are all its iterates.  Thus the $c$-values under consideration that may furnish a minimum value of $|\gamma|$ are $\pm 1$ and $-2f_0^8(0) \pm 1$, and plugging these into the expression for $\gamma$ gives 
\begin{equation} \label{gamma}
\gamma = \pm 2f_0^8(0) +  1 - m. 
\end{equation}
It is now easy to see that for $m \not\in \{-2,-1,0\}$, $|\gamma|$ is minimized by $m = 1$.  Indeed, the right-hand side of \eqref{gamma} is a polynomial in $m$; call it $g(m)$.  If its leading coefficient is positive, one checks that $g'(m) > 0$ for $m \geq 1$ and $g'(m) < 0$ for $m \leq - 3$ (one method is to use induction to examine the sign of $(f_0^8(0))'$).  The desired conclusion follows because $g(1)$ and $g(-3)$ are positive and $g(1) < g(-3)$.  A similar argument holds if the leading coefficient of $g(m)$ is negative.  Finally, having shown that $m = 1$, it follows that $|\gamma|$ is precisely the value given in \eqref{gam}.  


The condition (*) gives us more than just the fact that $\overline{f^8}$ is irreducible for some primes but $\overline{f^9}$ is not.  As in the introduction, the {\em arboreal Galois representation} attached to $f \in \Z[x]$ is the Galois group $G$ of the extension obtained by adjoining to $\Q$ all the preimages of $0$ under any iterate of $f$.  This set of preimages has a natural structure of a rooted tree, with the action of $f$ furnishing the connectivity relation.  The group $G$ has as quotient the Galois group $G_n$ of $f^n$ for any $n$, which acts naturally on the height-$n$ tree $T_n$ of preimages of $0$ under $f^n$.  By \cite[Theorem 3.3]{quaddiv}, condition (*) ensures that $G_8$ is as large as possible, i.e., the full tree automorphism group $\Aut(T_8)$.  This group contains elements acting on the roots of $f^8$ as a full $2^8$-cycle, which implies by the Chebotarev density theorem that there are primes for which $\overline{f^8}$ is irreducible.  On the other hand, the Galois group of $f^9$ is not as large as possible, since it contains no elements acting on the roots of $f^9$ as a $2^9$-cycle.  

\label{galoisdisc} This presents a contrast to the case of linear $\ell$-adic representations, i.e., Galois groups $G$ that are subgroups of $\GL_d(\Z_\ell)$.  Such representations arise from adjoining to the base field the coordinates of $\ell$-power torsion points on abelian varieties, or equivalently iterated preimages of the identity under multiplication by $\ell$.   The natural quotient giving the level-$n$ action is a subgroup of $\GL_d(\Z/\ell^n\Z)$.  In this case, if $G$ maps onto $\GL_d(\Z/\ell^n\Z)$ for certain small $n$, then $G$ must map onto $\GL_d(\Z/\ell^n\Z)$ for all $n$, and hence must be all of $\GL_d(\Z_\ell)$.  For instance, when $d = 2$ and $\ell \geq 5$, any $G \leq \GL_2(\Z_\ell)$ that surjects onto 
$\GL_2(\Z/\ell\Z)$ must be all of $\GL_2(\Z_\ell)$.  The salient difference is that the Frattini subgroup of $G \leq \GL_d(\Z_\ell)$ has finite index in $G$, while the Frattini subgroup of the automorphism group of the infinite tree of preimages of $0$ under a quadratic polynomial has infinite index.  For more on this, see \cite[Sections 3 and 5]{itend}.  For a more general discussion of surjectivity criteria for linear Galois representations, see \cite{vasiu}.  

To prove that such a phenomenon cannot occur in the present case of arboreal representations attached to quadratic polynomials, we would need to find, for each $n \geq 1$, some $f(x)$ satisfying condition (*) for $1 \leq i \leq n-1$ and also with $f^n(\gamma)$ a square.  While this remains out of reach, we are able to adapt the construction with $n = 9$ at the beginning of this section to show:

\begin{theorem} \label{primex}
Let $n \geq 2$.  Then there exists a quadratic $f \in \Z[x]$ that is stable, and such that $\overline{f^{n-1}}$ is irreducible for a positive density of primes, but $\overline{f^n}$ is reducible for all primes.   
\end{theorem}

Note that Theorem 4.1 immediately implies Theorem 1.2. The idea behind the proof of Theorem 4.1 is to choose $s = (f_0^{n-1}(0) - 1)^2$, rather than $s = (f_0^{n-1}(0) + 1)^2$ as was done in the construction at the beginning of this section.  The conclusion of Corollary \ref{qcor} still applies provided we can show that $-f(\gamma)$ is not a square, since $s$ is of the appropriate parity.  This choice gives 
\begin{align} 
\gamma & = (f_0^{n-1}(0) - 1)^2 - (f_0^{n-1}(0)^2 + m). \nonumber \\
& =  -2f_0^{n-1}(0) + 1 - m. \nonumber
\end{align}
and hence $-\gamma = 2f_0^{n-1}(0) + 1 > f_0^i(0)$ for $i = 1, \ldots, n-1$.  It follows that $f^i(\gamma) < 0$, and this allows us to circumvent having to verify condition (*), as the following lemma shows:


\begin{lemma} \label{lastlem}
Let $a_1, \ldots, a_k$ be negative integers, and let $q$ be a prime not dividing any $a_i$.  Then for any integer $c > 0$ with $q \nmid c$ there is a prime $p$ with $(qc/p) = -1$ and $(a_i/p) = -1$ for all $1 \leq i \leq k$, where $( \cdot / \cdot )$ denotes the Legendre symbol.  
\end{lemma}

\begin{remark}
Indeed, the set of $p$ with the desired property has positive density in the set of all primes.  
\end{remark}

\begin{proof}[Proof of Lemma \ref{lastlem}]
Let $r_1, \ldots, r_j$ be the primes dividing $|ca_1a_2 \cdots a_k|$, and note that by hypothesis none of the $r_i$ can equal $q$.  Using quadratic reciprocity, the Chinese remainder theorem, and Dirichlet's theorem on primes in arithmetic progressions, we may find a prime $p \equiv 3 \pmod{4}$ with $(r_i/p) = 1$ for all $1 \leq i \leq j$ and $(q/p) = -1$.  Then $(-1/p) = -1$, and it follows that $p$ is the desired prime.  
\end{proof}

If in fact the choice of $s = (f_0^{n-1}(0) - 1)^2$ caused each of $-f(\gamma), f^2(\gamma), f^3(\gamma), \ldots, f^{n-1}(\gamma)$ to be negative, then by Theorem \ref{fund} the prime $p$ in Lemma \ref{lastlem} would be the one required to prove Theorem \ref{primex}.  However, $-f(\gamma)$ is obviously positive in this case, and so we must do more.  
\begin{lemma} \label{reallast}
For each $n \geq 1$ there exist $m \in \Z$ and a prime $q$ with the following property.  Take $f_0(x) = x^2 + m$, 
$\gamma = -2f_0^{n-1}(0) + 1 - m$, and $f(x) = (x - \gamma)^2 + \gamma + m$.  Then $q$ divides $-f(\gamma)$ to the first power only and does not divide $f^i(\gamma)$ for any $i > 1$.  
\end{lemma}
Lemma \ref{reallast} is enough to establish Theorem \ref{primex}, since we may apply Lemma \ref{lastlem} with $c = -f(\gamma)/q$ and $a_i = f^{i+1}(\gamma)$ for $i = 1, \ldots, n-2$.  The resulting prime $p$ then has $\overline{f^{n-1}}$ irreducible, and by the proof of Corollary \ref{qcor} $f$ is stable and $\overline{f^n}$ is reducible for all primes.  

\begin{proof}[Proof of Lemma \ref{reallast}]
For $n = 1$ the statement is trivially true, so we begin with $n=2$.  Take $m = 3$. One checks directly that $f(\gamma) = -5$ and $f^i(\gamma) \equiv 4 \bmod{5}$ for all $i > 1$ (the latter can be done by calculating the orbit $f(\gamma), f^2(\gamma), \ldots$ modulo $5$).  Thus the lemma is true with $q = 5$.  

Suppose now that $n \geq 3$, and consider the case where $n \not\equiv 1 \bmod{3}$.  We claim that taking $m = 1$ and $q = 3$ suffices.  Note that with $m = 1$, the orbit $f_0(0), f_0^2(0), f_0^3(0), \ldots$ modulo $9$ is   
\begin{equation} \label{orbit}
1 \to 2 \to 5 \to 8 \to 2 \to \cdots
\end{equation}
Now $f(\gamma) = \gamma + m = -2f_0^{n-1}(0) + 1$, and so from \eqref{orbit} we have 
$f(\gamma) \equiv 3 \bmod{9}$ or $f(\gamma) \equiv 6 \bmod{9}$ for all $n \geq 3$ with $n \not\equiv 1 \bmod{3}$.  Thus $3$ divides $-f(\gamma)$ to the first power only.  Now observe that for any $i > 1$, $f^i(\gamma) - f(\gamma) = f_0^i(0) - m = (f_0^{i-1}(0))^2$, and hence any prime dividing both $f^i(\gamma)$ and $f(\gamma)$ must also divide $f_0^{i-1}(0)$.  It follows from \eqref{orbit} that $3 \nmid f^i(\gamma)$ for all $i > 1$.  

In the case where $n \equiv 1 \bmod{3}$, we take $m = 4$ and $q = 3$.  The orbit in \eqref{orbit} now becomes
$$4 \to 2 \to 8 \to 5 \to 2 \to \cdots,$$
and we argue as in the previous case.  
\end{proof}

\section{Results for function fields} \label{funcfields}

When $F$ is a function field over a finite field of odd characteristic, we cannot use the same proof as in Theorem \ref{numfield}, since now $2$ is a unit.  Indeed, there does not appear to be a stability result as general as that of Theorem \ref{numfield} that will allow us to mimic the construction of Corollary \ref{numfieldcor}.  However, it is still possible to give conditions on $m$ and $\gamma$ that ensure $f^n$ is irreducible but $\overline{f^n}$ is reducible for almost all primes.    

Let $F$ be a function field over a finite field $k$ of odd characteristic, and let $\calO$ be the integral closure of $k[t]$ in $F$.  In contrast with the usage of the previous two sections, we take a prime of $F$ to be slightly more general than simply the prime ideals lying in $\calO$.  Specifically, a prime of $F$ is a discrete valuation ring $R \subset F$ that contains $k$ and has field of fractions $F$.  Denote the maximal ideal of $R$ by $P$; we often refer to both $P$ and $R$ as a prime of $F$.  We may extend the valuation on $R$ to a multiplicative function $v_P : F^* \to \Z$, which we call the $P$-adic valuation.  
For all primes $P$ of $F$, the $P$-adic valuation satisfies the strong triangle inequality: for $x, y \in F^*$, $v_P(x + y) \geq \min\{v_P(x), v_P(y)\}$, with equality holding if $v_P(x) \neq v_P(y)$.

\begin{theorem} \label{funcfield}
Let $F$ be a function field over a finite field $k$ of odd characteristic, and let $\calO$ be the integral closure of $k[t]$ in $F$.  Let $m \in F$, and suppose that there are two primes $Q_1$ and $Q_2$ with $v_{Q_1}(m)$ positive, $v_{Q_2}(m)$ negative, and both odd.  Let $n \geq 3$, $f_0(x) = x^2 + m$, take $\gamma = m^{2^{n-1}} - f_0^n(0)$, and set $f(x) = (x- \gamma)^2 + \gamma + m$.  Then $f^n$ is irreducible over $F$ but $\overline{f^n}$ is reducible for each prime $P$ of $F$ with $v_P(m) \geq 0$.   
\end{theorem}

\begin{remark}
Note that $v_P(m) \geq 0$ for all but finitely many primes of $F$ \cite[Proposition 5.1]{Rosen}.  There also are (many) $m \in F$ that satisfy the hypotheses of Theorem \ref{funcfield}.  Indeed, fix a prime $Q_2$ of $F$ with $\deg Q_2$ odd, which is possible since $F$ has primes of all sufficiently large degrees by the Weil bound \cite[Theorem 5.12]{Rosen}.  Let $D_n$ be the divisor $nQ_2$.  For $n$ large enough, the Riemann-Roch theorem gives $l(D_n) - l(D_{n-1}) = 1$ \cite[p. 49]{Rosen}, where $l(D)$ is the dimension of the $k$-vector space $L(D) := \{m \in F^* : D + (m) \geq 0\} \cup \{0\}$.  Thus we may take $m \in L(D_n) \setminus L(D_{n-1})$ for $n$ odd and sufficiently large, whence $v_{Q_2}(m)$ is negative and odd.  Moreover, $Q_2$ is the only place with $v_{Q_2}(m) < 0$. By \cite[Proposition 5.1]{Rosen},
$$\sum_{\{P \, : \, v_P(m) < 0\}} -v_P(m) \deg{P} = \sum_{\{P \, : \, v_P(m) > 0\}} v_P(m) \deg{P},$$
and because the left-hand side is odd, the right-hand side is as well.  It follows that there must be a place $Q_1$ with $v_{Q_1}(m)$ positive and odd.  
\end{remark}

\begin{remark}
Unlike Theorems \ref{numfieldcor} and \ref{qcor}, the conclusion of Theorem \ref{funcfield} doesn't necessarily hold for $f^i$ with $i \geq n$.  Clearly $\overline{f^i}$ is reducible for any $i \geq n$ for each prime $P$ with $v_P(m) \geq 0$, but the lack of an equivalent of Theorem \ref{numfield} means we can't conclude that $f^i$ is irreducible over $F$.  Note that Proposition \ref{deddomstab} can't be used under the hypotheses of Theorem \ref{funcfield}, since 
$v_P(\gamma) = v_P(m)$ for all primes $P$ with $v_P(m) > 0$.    
\end{remark}

\begin{proof}[Proof of Theorem \ref{funcfield}]
Let $v_{Q_1}(m) = c_1  > 0$ with $c_1$ odd.  By the proof of Proposition \ref{deddomstab}, $v_{Q_1}(f_0^i(0)) = c_1$ for all $i \geq 1$, and hence 
$$v_{Q_1} \left(\frac{f^{n-1}(\gamma) -\gamma \pm \sqrt{f^n(\gamma)}}{2}\right) = v_{Q_1}(f_0^{n-1}(0) \pm m^{2^{n-2}}) = c_1,$$
where the last equality follows from the strong triangle inequality and the assumption that $n \geq 3$.  
Hence neither of $(f^{n-1}(\gamma) -\gamma \pm \sqrt{f^n(\gamma)})/2$ is a square in $F$.  

Let $v_{Q_2}(m) = c_2 < 0$ with $c_2$ odd.  Note that $f_0^n(0) = m^{2^{n-1}} + 2^{n-2}m^{2^{n-1} - 1} + \cdots$, and thus $v_{Q_2}(\gamma) = (2^{n-1} - 1)c_2$, which is odd.  Moreover, for $i < n$, 
$$v_{Q_2}(f_0^i(0)) = (2^{i-1})v_{Q_2}(m) >  v_{Q_2}(\gamma),$$  
where the final inequality follows because $n \geq 3$ ensures $2^{i-1} < 2^{n-1} - 1$.
Because $f^i(\gamma) = f_0^i(0) + \gamma$, it follows that 
$v_{Q_2}(f^i(\gamma)) = v_{Q_2}(\gamma)$, and hence $f^i(\gamma)$ is not a square in $F$.   Hence by Theorem \ref{altfund}, $f^n$ is irreducible over $F$.  On the other hand, if $P$ is a prime of $F$ with $v_P(m) \geq 0$, then also $v_P(m^{2^{n-1}}) \geq 0$.  The proof of Theorem \ref{deddom} then shows that $\overline{f^n}$ is reducible.  
\end{proof}

When $F = k(t)$, we can simplify the hypotheses of Theorem \ref{funcfield}.  Recall that in this case there is a prime $P_\infty$ given by the discrete valuation ring $k[t^{-1}]$, whose unique maximal ideal is generated by $t^{-1}$.  The corresponding valuation $v_{P_{\infty}}$ attaches the value $\deg(g) - \deg(f)$ to the element $f/g \in F$.  We refer to a prime $P$ of $F$ with $P \neq P_\infty$ as a finite prime.  

\begin{corollary} \label{ratffcor}
Let $k$ be a finite field of odd characteristic, $F = k(t)$, $\calO = k[t]$, and suppose that $m = f(t)/g(t) \in F$ with $(f,g) = 1$, $\deg(f)$ odd, $\deg(g)$ even, and $\deg(f) > \deg(g)$.  Let $n \geq 3$, $f_0(x) = x^2 + m$, take $\gamma = m^{2^{n-1}} - f_0^n(0)$, and set $f(x) = (x- \gamma)^2 + \gamma + m$.  Then $f^n$ is irreducible over $F$ but $\overline{f^n}$ is reducible for each finite prime $P$ of $F$ with $v_P(g) = 0$.
\end{corollary}

\begin{proof}
By hypothesis $v_{P_{\infty}}(m) = \deg(g) - \deg(f)$ is negative and odd.  Because $f$ has odd degree, it cannot be a constant times a square, and hence there is a prime $P$ with $v_P(f)$ positive and odd.  But $(f,g) = 1$, and thus $v_P(f) = v_P(m)$, and the hypotheses of Theorem \ref{funcfield} are satisfied. 
\end{proof}

To illustrate Corollary \ref{ratffcor}, let $n = 3$ and $m = t$.  Then $\gamma =  t^4 - (t^4 + 2t^3 + t^2 + t) = -2t^3 - t^2 - t$.  Take 
$$f(x) = (x - \gamma)^2 + \gamma + t = x^2 + (4t^3 + 2t^2 + 2t)x + 4t^6 + 4t^5 + 5t^4.$$
Then $f^3(x)$ is irreducible over $F$ but reducible modulo all finite primes of $F$.  In other words, for any $c$ in the algebraic closure of $k$, the specialization of $f^3(x)$ at $t = c$ is reducible over $k(c)$, even though $f^3(x)$ is irreducible over $F$.  

We note that Theorem \ref{funcfield} doesn't apply when $n = 2$, since then $f_0^n(0) = m^2 + m$, which means according to the recipe of Theorem \ref{funcfield}, $\gamma = -m$.  But then $\gamma + m = 0$, and so $f$ is reducible.  However, this may be remedied by choosing $r$ with $r/2$ a non-quadratic residue in $k$ and taking $\gamma = (m + r)^2 - m^2 - m$.  Then $f^2(\gamma) = \gamma + m^2 + m = (m+r)^2$.  Moreover, $-f(\gamma) = -(\gamma + m) = -(2rm + r^2)$.  Because $r/2$ is not a quadratic residue, $r \neq 0$, and thus $-(2rm + r^2)$ has odd $Q_2$-adic valuation (under the hypotheses of Theorem \ref{funcfield}), and so is not a square in $F$.  Therefore $f$ is irreducible.  Finally, we have 
$$\frac{-m + \sqrt{f^2(\gamma)}}{2} = \frac{r}{2},$$
which is not a square in $F$, showing that $f^2$ is irreducible by Theorem \ref{altfund}.  It is worth noting that if we extend the field of constants of $F$ to be $k(\sqrt{r/2})$ then $f^2$ becomes reducible.  


\section{the number of stable primes}  \label{numstab}

The purpose of this section is to investigate, for given monic, quadratic $f$ defined over a global field $F$, the number of primes of $F$ for which $\overline{f}$ is stable.  For simplicity let us suppose that $f$ is defined over $\calO$, which we take to be the ring of integers of $F$ in the number field case and the integral closure of $k[t]$ in the case where $F$ is a function field over the finite field $k$ (of odd characteristic).  Then $f(x)$ may be written as $(x - \gamma)^2 + \gamma + m$, with $\gamma \in \frac{1}{2}\calO$ and $m \in \frac{1}{4}\calO$.  In the function field case the reductions $\overline{\gamma}$ and $\overline{m}$ are defined for all primes not lying over $P_\infty$, while in the number field case they are defined for all primes not lying over $2$.  For the latter, $\overline{f}$ cannot be stable, as indeed its third iterate must always be reducible \cite{ahmadi}.  

Recall that the \textit{affine span} of a subset $S$ of a vector space $V$ is the collection of all $v \in V$ that can be written as a linear combination of elements of $S$ whose weights sum to 1.  Of interest here is the $\Z/2\Z$-vector space $F^*/F^{*2}$. If $S = \{s_1, s_2, \ldots\} \subseteq F^*$, then the affine span of $S$ (considered as a subset of $F^*/F^{*2}$) is the collection of all $F^{*2}$-cosets with a representative of the form $\prod_{j \in J} s_j$,
where the number of elements in the set $J$ is odd.  Note that the affine span of $S$ contains the origin (i.e., the identity coset) if and only if a product of an odd number of elements of $S$ is a square in $F$.

\begin{theorem} \label{stabnum}
Let $F$ be a global field, and $f \in F[x]$ monic and quadratic with critical point $\gamma$.  Let 
$S = \{-f(\gamma), f^2(\gamma), f^3(\gamma), \ldots\}$.  
\begin{enumerate}

\item If $0 \in S$ or if $0 \not\in S$ and the affine span of $S$ in $F^*/F^{*2}$ contains the origin, then there is an iterate of $f$ that is reducible modulo all primes.  
\item If $ 0 \not\in S$ and the affine span of $S$ in $F^*/F^{*2}$ is finite, say of cardinality $2^d$, and does not contain the origin, then $\overline{f}$ is stable for a set of primes of density $2^{-d-1}$.  
\end{enumerate}
\end{theorem}

Note that in Theorem \ref{stabnum}, (1) implies that $\overline{f}$ is stable for no primes, while (2) implies $\overline{f}$ is stable for infinitely many primes.   In assertion (2), we use the notion of natural density for sets of primes in number fields and Dirichlet density for sets of primes in function fields.  
In the case $F = \Q$, the positive-density set of primes referenced in (2) is by quadratic reciprocity the union of congruence classes for some fixed modulus.  

\begin{example}
Let $m = -1$ and $\gamma = -1$, so that $f(x) = (x+1)^2 - 2$.  Then $S = \{2, -1, -2\}$, and we have the relation $s_3 = s_1s_2$.  Multiplying through by $s_3$ makes clear that the affine span of $S$ contains the origin.  Note that $f^3$ is reducible modulo all primes.  
\end{example}

\begin{proof}[Proof of Theorem \ref{stabnum}]
Suppose first that $0 \in S$.  Then $\gamma$ is a root of $f^n(x)$ for some $n \geq 1$, and thus 
$(x- \gamma) \mid f^n(x)$, so that $f^n(x)$ is reducible in $F[x]$. Therefore $f^n(x)$ is reducible modulo all primes. For the remainder of assertion (1), choose $n$ large enough so the first $n$ elements $s_1, \ldots, s_n$ of $S$ satisfy some equality 
\begin{equation} \label{subprod}
r^2 = \prod_{j \in J} s_j
\end{equation}
and $\#J$ is odd.  Then there can be no $\p$ with all $s \in S$ non-squares modulo $\p$, since then $\prod_{j \in J} s_j$ would be a non-square modulo $\p$, which is absurd.  By Theorem \ref{fund} we thus have $f^n$ reducible modulo all primes.  

For assertion (2), note that by Theorem \ref{fund}, the set of primes $\p$ such that $\overline{f}$ is stable for $\p$ coincides with the set $T$ of primes $\p$ such that no element of $S$ is a square in $\calO/\p$.  Consider the extension $E$ of $F$ obtained by adjoining to $F$ the square roots of all elements of $S$.  Then $E$ is a finite Galois extension of $F$ with $\Gal(E/F)$ an elementary abelian $2$-group.  Moreover, $\p \in T$ if and only if $\Frob_\p \in \Gal(E/F)$ maps $\sqrt{s}$ to $-\sqrt{s}$ for each $s \in S$.  

By Kummer theory, $|\Gal(E/F)|$ is the size of the span of $S$ in the vector space $F^*/F^{*2}$.  Let $B \subseteq S$ be a basis for ${\rm Span}(S)$. Each $s \in S \setminus B$ must be a square times the product of an odd number of elements of $B$, for otherwise multiplying both sides by $s$ gives an equality as in \eqref{subprod}, with $\#J$ odd. This contradicts our supposition that the affine span of $S$ does not contain the origin.

It follows now that the affine span of $S$ consists of the $F^{*2}$-cosets whose representatives are products of an odd number of elements of $B$. Thus the affine span of $S$ has half as many elements as the span of $S$, and hence we have $\#{\rm Span}(S) = 2^{d+1}$.  Moreover, the observation that each $s \in S \setminus B$ must be a square times the product of an odd number of elements of $B$ implies that the unique $\sigma \in \Gal(E/F)$ with $\sigma(\sqrt{b}) = -\sqrt{b}$ for all $b \in B$ also satisfies $\sigma(\sqrt{s}) = -\sqrt{s}$ for all $s \in S$. 
By the Chebotarev density theorem (see \cite[p. 545]{neukirch} for the number field case, \cite[p. 125]{Rosen} for the function field case), the density of $\p$ with $\Frob_\p = \sigma$ is $2^{-d-1}$. 
\end{proof}

\begin{conjecture} \label{stabconj}
Let $F$ be a global field, and $f \in F[x]$ monic and quadratic with critical point $\gamma$.  Let 
$S = \{-f(\gamma), f^2(\gamma), f^3(\gamma), \ldots\}$.  If $0 \not\in S$ and the affine span of $S$ as a subset of $F^*/F^{*2}$ is infinite and does not contain the origin, then $\overline{f}$ is stable for only finitely many primes.
\end{conjecture}

Note that under the hypotheses of Conjecture \ref{stabconj}, it follows from Kummer theory and the Chebotarev density theorem that the set of $\p$ for which $\overline{f}$ is stable has density zero. Conjecture \ref{stabconj} appears difficult to prove.  However, the following heuristic suggests that it is true.  For $\p \in \calO$, denote by $N_\p$ the the number of elements of $\calO/\p:=F_\p$.  We need two main assumptions: that the elements of the orbit of $\overline{\gamma}$ behave like a random orbit of a random self-map of $F_\p$ and that the elements of $S$ are multiplicatively independent.  The orbit of a random point under a random self-map of $F_\p$ has length bounded below by $\sqrt{N_\p}$ \cite{Flajolet} (see also \cite[Section 6]{jhsnyjm}).  Hence $\overline{f}$ is stable for $\p$ if none of $-f(\gamma), f^2(\gamma), f^3(\gamma), \ldots, f^j(\gamma)$ is a square in $F_\p$, for some $j \geq \sqrt{N_\p}$.  As in the proof of Theorem \ref{stabnum}, part (2), the set of primes for which this is true has density $1/r$, where 
$r$ is the size of the span of $-f(\gamma), f^2(\gamma), f^3(\gamma), \ldots, f^j(\gamma)$ in $F^*/F^{*2}$.  By our independence assumption, $r = 2^j$, and so the ``probability" that $\overline{f}$ is stable is at most $2^{-\sqrt{N_\p}}$.  This gives that the expected number of primes for which $\overline{f}$ is stable is
\begin{equation} \label{norm2}
\sum_\p 2^{-\sqrt{N_\p}}.
\end{equation}
When $F$ is a number field, let $d = [F:\Q]$, and note that for a given rational prime $p$, the sum \eqref{norm2} taken over $\p \mid (p)$ can be at most $d/2^{\sqrt{p}}$, which occurs when $(p)$ splits completely in $F$.  Hence the full sum in \eqref{norm2} is at most $\sum_p d/2^{\sqrt{p}}$, which is less than $d \sum_n 1/2^{\sqrt{n}}$.  Separating this last sum into the pieces $i^2 \leq n \leq (i+1)^2 - 1$, we see that it is bounded above by $d \sum_i (2i+1)/2^i$, which converges.  A similar argument holds in the function field case.

It would be very interesting to establish the conclusion of Conjecture \ref{stabconj} for any single polynomial.  We consider the case of $F = \Q$, $f(x) = x^2 + 1$.  Odoni  \cite{odoni} first observed that $\overline{f}$ is stable for $p = 3$, and also remarked on the central role that the sequence $-f(\gamma), f^2(\gamma), f^3(\gamma), \ldots$ plays in the Galois theory of iterates of $f(x)$.  His work paved the way for Stoll's proof  \cite{stoll} that the arboreal representation attached to $f(x)$ is surjective, i.e., the Galois groups of iterates of $f(x)$ are as large as possible.  
\begin{conjecture} \label{concrete}
Let $F = \Q$ and $f(x) = x^2 + 1$.  Then $\overline{f}$ is stable for $p = 3$ and for no other primes.
\end{conjecture}

Note that for $f(x) = x^2 + 1$, the set $\{-f(0), f^2(0), f^3(0), \ldots\}$ is linearly independent over $\Q^*/\Q^{*2}$ \cite{stoll}, and in particular its affine span is infinite and does not contain the origin.  

Using a computer algebra system such as MAGMA, one computes that the first 20 elements of $-f(\gamma), f^2(\gamma), f^3(\gamma), \ldots$ are all non-squares modulo $p$ for 42 of the $50,847,534$ primes $\leq 10^9$.  Apart from 3, each of these primes has $f^n(\gamma)$ a square modulo $p$ for some $n \leq 25$, thereby verifying Conjecture \ref{concrete} for primes $\leq 10^9$.  As further evidence, we give the following result, though we first define some terminology.  Let $a, f(a), f^2(a), \ldots$ be a finite orbit, and take $f^0(a) = a$. Let $r$ be the minimal positive integer with $f^r(a) = f^s(a)$ for some $0 \leq s < r$.  Then the \textit{tail} of the orbit is $a, f(a), \ldots, f^{s-1}(a)$ when $s > 0$, and is empty otherwise.  By the length of the tail, we mean $s$.   
\begin{proposition}
Let $f(x) = x^2 + 1$, and suppose that $\overline{f}$ is stable for a prime $p$.  Then the orbit of $0$ under $\overline{f}$ has tail of length two.  
\end{proposition}  

\begin{proof}
To ease notation, let $a_n = \overline{f}^n(0)$ for $n \geq 0$, and note that $a_0 = 0$ and $a_n = a_{n-1}^2 + 1$ for $n \geq 1$.  Let $r$ be minimal with $a_r = a_s$ for some $s < r$.  If $s = 0$ then 
$a_r = 0$, and hence $\overline{f}$ is not stable.  If $s = 1$ then $a_r = 1$, and hence $a_{r-1} = 0$, so again $\overline{f}$ is not stable.  So assume $s \geq 2$.  
Then $a_{r-1}^2 = a_{s-1}^2$.  But by the minimality of $r$, we must have $a_{r-1} \neq a_{s-1}$.  Hence $a_{r-1} =  -a_{s-1}$. 
Note that not all of $a_{s-1}, -a_{s-1},$ and $-1$ can be non-squares in $\Z/p\Z$.  Because $-1 = -a_1$, this shows that $a_{s-1}, a_{r-1},$ or $-a_1$ is a square in $\Z/p\Z$.  If the square is $a_{r-1}$ or $-a_1$, or if $s > 2$, then one of $-a_1, a_2, a_3, \ldots$ is a square, and $\overline{f}$ is not stable by Theorem \ref{fund}.  Therefore if $\overline{f}$ is stable then $s = 2$.  \end{proof}

We note that if $s = 2$, then $\overline{f}^r(0) = 2$ for some $r \geq 2$, and indeed $\overline{f}^{r-1}(0) = -1$, since otherwise $\overline{f}^{r-1}(0) = 1 = \overline{f}^1(0)$, contradicting $s = 2$.  Thus the only $p$ for which $\overline{f}$ has a chance of being stable are those with $f^n(0) \equiv -1 \bmod{p}$ for some $n$.  By factoring $f^n(0) + 1$ for $1 \leq n \leq 9$, one sees that apart from $3$, all primes with $f^n(0) \equiv -1 \bmod{p}$ for $1 \leq n \leq 9$ are congruent to $1$ modulo 4, and thus $-1$ is a square modulo $p$, so already $\overline{f}(x)$ is reducible.  However, there are factors of $f^n(0) + 1$ with $n = 10, 11$ that are congruent to $3$ modulo 4.  

We note that $f^n(0) + 1$ may be obtained from $f^{n-1}(0) + 1$ by applying $g(x) = (x-1)^2 + 2$.  Indeed, $g^{n}(1) = f^{n}(0) + 1$, so the only primes for which $x^2 + 1$ has a chance of being stable are those dividing some element of the forward orbit of the critical point of $g$.   

As a final remark, we note that for general monic, quadratic $f \in F[x]$, there is presently no good method for determining the infinitude of the affine span of the set $S$ in Conjecture \ref{stabconj}. One exception is in the cases where $f^n(\gamma)$ is a rigid divisibility sequence or the orbit of $0$ under $f$ is finite. In these cases one can prove that for infinitely many $n$, there is a prime dividing $f^n(\gamma)$ with odd multiplicity but not dividing $f^i(\gamma)$ for any $i < n$ (see \cite{quaddiv} for details). This implies that the affine span of $S$ is infinite.

\section*{acknowledgements}
The author thanks the anonymous referee for valuable comments and suggestions, including improved statements of many of the results and conjectures of Section \ref{numstab}.

\bibliographystyle{plain}

\end{document}